\documentclass[A4paper,11pt,authoryear]{elsarticle}

\usepackage{amssymb}
\usepackage{amsmath}
\usepackage{url}
\usepackage[hidelinks]{hyperref}
\usepackage{booktabs,array,makecell,adjustbox,amsmath,bm,subcaption}

\newtheorem{definition}{Definition}

\newcommand{\ci}[2]{\footnotesize[\,#1; #2\,]}
\newcolumntype{L}[1]{>{\raggedright\arraybackslash}p{#1}}
\newcolumntype{C}[1]{>{\centering\arraybackslash}p{#1}}
\newcolumntype{R}[1]{>{\raggedleft\arraybackslash}p{#1}}

\begin{document}

\begin{frontmatter}

\title{Maximizing Effectiveness and Equity in Kidney Exchange Programs for Novel Compatibility Paradigms}
\author{ Valentina Peralta Clarke, Faculty of Engineering, Universidad Adolfo Ib\'a\~nez, Santiago de Chile, Chile \\ 
Hans de Ferrante, Eurotransplant, Leiden, The Netherlands \\
Francisco Perez Galarce, Universidad de Las Americas, Faculty of Engineering and Business, Santiago de Chile, Chile \\
Joris van de Klundert*, Business School Universidad Adolfo Ib\'a\~nez, Vina del Mar, Chile \\
* corresponding author, email: joris.vandeklundert@uai.cl}

\begin{abstract}
Kidney Exchange Programs (KEPs) promote access to living donor transplantation for patients suffering from end stage renal disease. The HLA compatibility between transplant recipients and donors plays an important role when solving the matching problems occuring in KEPs. Compatibility defines the feasible solution space and often occurs in a weighted form in the objective function. Recently, the paradigms used to express HLA compatibility have advanced substantially from antigens to alleles and eplets. 

We show how KEP effectiveness and equity vary with these three paradigms. The equity analysis focuses on ethnic inequities and the result confirm that moving from antigen to allele and eplet paradigms may enlarge inequities. We present new optimization models based on equity weighting that remedy ethnic inequities for all three paradigms without sacrificing access to transplant.


\end{abstract}

\end{frontmatter}

\newpage

\section{Introduction}
\label{sec: introduction}

Living donor kidney transplantation (LTx) is the most cost-effective and therefore preferred treatment for patients suffering from end stage reneal disease (ESRD) \citep{Wolfeetal2009, Mohnenetal2019}. The resulting transplant outcomes are strongly associated with the patient and donor being blood type (ABO) compatible and Human Leucocyte Antigen (HLA) compatible \citep{takemoto2004hla}. Incompatibility may cause a rejection response that can cause graft failure and, therefore, return to dialyis or a next transplant. The immunological understanding of the relationship between ABO and HLA compatibility on the one hand and rejection on the other hand has advanced considerably over the past decades \citep{mattoo2024improving, tambur2021significance}. As the nuances of this relationship are better understood, compatibility is increasingly defined along a continuum ranging from fully compatible to completely incompatible, instead of dichotomously as compatible versus incompatible \citep{mattoo2024improving}. 

If a patient has found a willing living donor who is perceived to be insufficiently compatible, the patient-donor pair may join a Kidney Exchange Program (KEP). KEPs facilitate the exchange of donors to improve patient donor compatibility and therefore LTx effectiveness.

Compatibility thus plays an essential role in the optimization of KEPs. Donor patient combinations classified as insufficiently compatible are eliminated from the {\it compatibility graph} that defines the feasible solution space. Moreover, the continuously defined compatibility among the remaining feasible combinations may be expressed in the weights of the objective function of optimization models to solve the matching problems arising in KEPs \citep{biro2021modelling}. In this research, we analyse how recent advances in the understanding of HLA compatibility impact the effectiveness and equity of KEPs and how inequities can be addressed.\\




The technologies and theory to capture HLA compatibility are advancing rapidly \citep{huang2019assessing}. Historically, HLA compatibility has been defined as the number of mismatches at the antigen level for the HLA-A, -B, and -DR loci. At present, however, the HLA typing of donors and patients is often done for up to eleven antigens at the HLA-A, -B, -C, -DR, -DQ, and -DP loci \citep{tambur2023can, jaramillo2023human}. The relative importance assigned to the individual HLA loci has evolved, with HLA-DR and -DQ increasingly considered as important \citep{tambur2023can}. 

Compatibility definitions based on antigen level data (also known as {\it low resolution} data) are commonly adopted and associated with transplant effectiveness. However, the prediction accuracy of low resolution based survival prediction models is modest \citep{van2025comparative}. The use of low resolution data has also been associated with inequity in access to transplantation, disadvantaging patients with highly sensitized immune systems and ethnic minorities \citep{takemoto2004hla, karahan2025equitable}. Ethnic inequity has been a long-term concern in the US and is increasingly important elsewhere, including Europe, as migrant populations grow. 

Compared to antigen level data, the more refined  {\it allele level} data (or {\it high resolution}  data) can significantly better explain transplant failure \citep{jaramillo2023human} (see \ref{subsec:HLAmatching}). Allele level data are presently routinely considered in a variety of allocation systems across the globe \citep{mankowski2024balancing, tambur2023can,  karahan2025equitable}. 

A more recent development is to explain immune responses based on {\it epitopes}. Epitopes are amino acid structures occurring at the aforementioned loci and recognized by the immune system \citep{duquesnoy2001hlammatchmaker} (see \ref{subsec:HLAmatching}). Functional epitopes, known as {\it eplets} can significantly better predict graft failure than antigen level data \citep{wiebe2019hla} . Eplets can improve transplant effectiveness for patients whose antigen or allele level mismatch is modest while the eplet mismatch is high \citep{sapir2020epitopes, mattoo2024improving}. 


These advancements in compatibility paradigms feed into the debate on the survival improvements they facilitate and of the consequences of using survival predictions in optimization methods for kidney allocation such as in KEPs \citep{tambur2023can, sapir2020epitopes, mankowski2024balancing}. With two small-scale case studies, present evidence of the effects of applying eplet-based allocation mechanisms in KEPs is very scarce \citep{kausman2016application,tafulo2020improving}. The KEP of the US based National Kidney Registry (NKR) uses eplet-based matching \citep{NKR2025} but details and results have not appeared in the scientific literature yet. \cite{lemieux2021matchmaker} argue that the understanding of eplet-based matching is still in flux and that it is unclear if and how eplet-based matching might improve KEP effectiveness. Our first research objective therefore is to assess how KEP effectiveness varies with the adoption of the three aforementioned compatibility paradigms. 

The effectiveness improvement attainable by using more advanced compatibility paradigms in organ allocation programs such as KEPs might outweigh possible adverse equity impacts \citep{kausman2016application, mangiola2024hla, mankowski2024balancing}. On the other hand, the scarce evidence causes the equity impact to be insufficiently understood and calls for caution regarding the implementation of eplet-based compatibility in the optimization of kidney allocation systems such as KEPs \citep{tambur2023can}. The second research question therefore seeks to identify the equity impact on access to transplant and quality of transplant in KEPs when moving from antigen- to allele- and eplet-based matching. As the distributions of antigens, alleles and eplets vary in different ways among ethnic groups, the focus is especially  on ethnic equities in allocation. The third research question studies how to resolve ethnic inequities in access to transplant and quality of transplants associated with the three HLA compatibility paradigms. 

Equity has been previously considered in the operations research literature on kidney allocation and KEP optimization\citep{bertsimas2013fairness, dickerson2014price, van2022eliminating, st2025adaptation}. The operations research literature has so far, however, not contributed to addressing differences in effectiveness and equity associated with the adoption of more advanced compatibility paradigms. Our research thus presents optimization methods and a first case study on the adoption of advanced compatibility paradigms in the optimization of a large KEP.

\section{Methods} \label{sec:methods}

\subsection{Effectiveness and Equity of Transplantation} \label{subsec:accessandoutcomemeasures}

Reports of leading organ allocation programs present performance measures such as transplant probability and post transplant survival over one or multiple years for participating patients as well as for various subpopulations, e.g. per gender, age group, ABO blood type, et cetera \citep{nts_annual_2021, lentine2024optn, uk_annual_20232024}. This subsection introduces these measures and corresponding notation as a basis to define effectiveness and equity. 

Let $\mathcal{T} = [T_0,T_E]$ be a time period for which KEP performance is considered. Let $\mathcal{P} = {p_1,\ldots,p_m}$ be the set of all pairs of patients and donors participating in the KEP during $\mathcal{T}$. Focusing on the transplantation intention of the KEP, the remainder uses the term recipient rather than patient regardless of whether a transplant is received. Each pair $p_i, i = 1\ldots,m$, is defined by a quadruple $(r_i,d_i,a_i,l_i)$ where $r_i$ is the recipient (patient), $d_i$ is the donor, $a_i$ is the time of enrollment, and $l_i$ is the time pair $p_i$ leaves if not matched before. Reasons to leave a KEP before being matched are various and include poor health, death, and receiving a transplant outside of the KEP \citep{nts_annual_2021}. 

In the remainder, $T_0 < l_i$ and $a_i < T_E$, for all $i=1,\ldots,m$. Thus, some pairs may already have joined the KEP at $T_0$ and some pairs may still be in the KEP at $T_E$. Moreover, we disregard practical challenges that might cause matches proposed in solutions obtained after KEP optimization not to result in actual transplantation. Therefore the (maximum) total number of recipients $M$ matched by time $T_E$ equals the (maximum) total number of transplants in the KEP over time period $\mathcal{T}$. 


With the view that kidney transplant is life saving, the total number of transplants is the most commonly reported performance measure in research and practice \citep{smeuldersetal2021}. This measure disregards the distribution of the matching over relevant subpopulations and the health effects of transplantation as they may vary with the HLA compatibility among the donors with whom recipients can be matched. The total number of transplants can be viewed to be a process measure rather than an outcome measure that expresses KEP effectiveness \citep{cookson2020distributional, smeuldersetal2021, van2025effectiveness}. 

LYFT is an outcome measure to express KEP effectiveness. It is defined as the total number of life years from transplant (LYFT) for recipients who receive a transplant  relative to their life years without transplant \citep{wolfe2008calculating, segev2009evaluating}. As the counterfactual remains unknown, LYFT cannot be measured and has to be estimated. The same applies to the more advanced expected health utility gained as measured by expected quality adjusted life years gained (QALYs). 

Expected LYFT and QALYS gained can be estimated by survival models. While post kidney transplant survival models have been widely studied and inform allocation policies, their modest accuracy has been an argument to call for caution with their practical use \citep{segev2009evaluating, tambur2023can}. Moreover, (long term) survival data after living donor kidney transplant in relation to allele and eplet level HLA compatibility still are rather scarce \citep{tambur2023can}. Transplant survival models that use these refined and developing paradigms are therefore in early stages of development. Until more accurate HLA compatibility based survival models are available for relevant subpopulations, process measures that capture HLA compatibility at time of transplant can serve as a proxy for KEP effectiveness and equity and we follow this approach in the remainder \citep{smeuldersetal2021, uk_annual_20232024, lentine2024optn}.


For recipients who receive a transplant, waiting time until transplant forms another important process measure. Other recipients may still be in the KEP or may have left. To address the latter, process measures such as the probability to receive a transplant at all and the percentage of pairs that leave before receiving a transplant in the KEP often are reported as well. \\


Let $S = \{s_1,\ldots,s_{|S|} \}$ be a partitioning of the KEP population $p_1,\ldots,p_m$ into subpopulations. Then, the following KEP performance measures summarize and operationalize the above:
\begin{enumerate}
    \item Fraction $F(s)$ 
    of the members of subpopulation $s \in S$ that are matched during $\mathcal{T}$, 
    \item Fraction $L(s)$  
    of the members of subpopulation $s \in S$ that leave the KEP during $\mathcal{T}$ without being matched, 
    \item The fraction $1 - F(s) - L(s)$ of patients of a subpopulation $s \in S$ remaining in the KEP at time $T_E$. This measure is relevant as accumulation of (hard-to-match) patient populations in the KEP pool is a persistent KEP performance issue.
    \item Average compatibility score $HLA(s)$ achieved for the members of subpopulation $s \in S$ that have been matched during $\mathcal{T}$, 
    \item Average waiting time $W(s)$ achieved for the recipients in subpopulation $s \in S$ that have been matched $\mathcal{T}$. The waiting time is defined for matched pairs as the time elapsed between the arrival time $a_i$ and the time at which $p_i$ is matched (see below).
\end{enumerate}
When reported for the entire population $\mathcal{P}$, these measures can serve to define KEP effectiveness. When considered for proper subpopulations $s \subset P$, such as age groups, according to gender, ethnic groups, ABO blood types, or combinations of these demographic dimensions, these measures can serve as a basis to define equity, as further elaborated in subsection \ref{subsec:EquityWeighting}.

\subsection{HLA Compatibility} \label{subsec:HLAmatching}

HLA (in)compatibility comes in two basic forms. First, the immune system of a recipient may already have learned to recognize some specific HLA antigens, alleles, and eplets as non-self, because of earlier exposure to foreign blood or tissue. The immune system might then be capable of producing donor-specific antibodies (DSAs) that result in an immediate response to reject the transplanted organ. Donor antigens, alleles, or eplets which the recipient immune system has already learned to recognize as non-self and for which it therefore has DSAs are often perceived to be forbidden and as a reason to classify a donor as incompatible. The strength of such rejection responses may vary, however, and standards for the appraisal of DSA data are only starting to emerge \citep{mattoo2024improving}. 
This research considers DSA-based incompatibility as given. 

Second, for each of the HLA loci, there may be a difference between the antigens, alleles, or epitopes of the recipient and the donor. The recipient's immune system may learn to recognize these differences over time (months or even years) and develop \textit{de novo} donor-specific antibodies (dnDSAs). HLA mismatch refers to a difference in HLA for one or more loci, whether at the antigen, allele, or eplet level. For each of these three levels, the probability of kidney failure as a result of graft rejection increases with the mismatch grade \citep{takemoto2004hla, mattoo2024improving}. However, the strength of the relationship between mismatch grade and graft survival differs across the loci. Current evidence emphasizes the importance of the HLA-B, HLA-DR and HLA-DQ loci \citep{tambur2021significance}.

Traditionally, mismatches have been determined at the antigen level, i.e., using low-resolution data. Currently, 128 antigens are known across the A, B, C, DR, and DQ loci \citep{jaramillo2023human, barker2023ipd}. 

Allele level data, also know as high resolution data, are more refined than the low level antigen data and can partly explain the variance in the relationship between antigen level data and graft survival.  For a recipient who has DSAs against a donor's HLA when considering antigen level data, it may therefore turn out that the DSAs may not recognize the donors alleles and that the antigen level mismatch vanishes at the allele level \citep{mattoo2024improving}. Conversely, a recipient may develop  dnDSAs against an allele present in the donor, which was not mismatched at the antigen level \citep{Tambur2018}. Hence, the view that high resolution typing is more adequate to assess compatibility for kidney transplantation than low resolution typing \citep{jaramillo2023human}. HLA typing is increasingly common in practice and more than 30.000 alleles are presently known for the A, B, C, DR and DQ loci \citep{jaramillo2023human, barker2023ipd}. 

\begin{figure}
    \centering
    \includegraphics[width=1.0\linewidth]{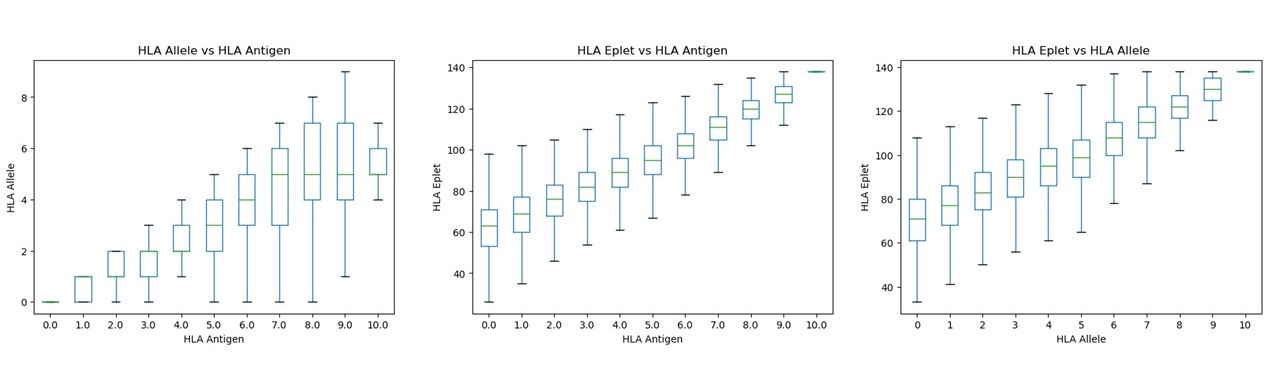}
    \caption{Allele level match per antigen level match for the case study population}
    \label{fig:paragimcorrelations}
\end{figure}

Low resolution (antigen) level can be accurately derived from the high resolution (allele) level typing as is applied in the case study reported below. The opposite is not true and probabilistic imputation of antigen level data to allele level data is relatively inaccurate for ethnically diverse populations and therefore should be used with caution \citep{engen2021substituting}. Figure \ref{fig:paragimcorrelations} shows the relationship between the number of allele level matches and number of antigen level matches for the population of the case study (Section \ref{sec:casestudy}). 

Epitopes are defined as the part of HLA molecules to which antibodies bind. Eplets are patches of surface-exposed amino acids of the HLA-protein, and are surrogates for these epitopes. These eplets are derived from computational studies, and are catalogued in curated databases. Eplets may be unique to a single antigen or allele or may be shared across antigens at different loci. High-resolution data are required to determine which eplets are present in the donor and recipient. The {\it eplet mismatch load} can then be quantified by comparing these eplets \citep{duquesnoy2001hlammatchmaker}. While the immunological validation of eplets is still a work in progress \citep{tambur2023can}, evidence shows eplets to enable more accurate prediction of graft failure compared to antigen or allele level based predictions \citep{wiebe2013class, tafulo2021hla, mattoo2024improving}. 

With a maximum of 67 for the DR and DQ loci and 138 when also considering the HLA A,B and C loci, the eplet mismatch load reported in the case study below is much higher than the maximum number of mismatched loci. The relationships of eplet mismatch load with antigen and allele level HLA match for the case study data are also depicted in Figure \ref{fig:paragimcorrelations} (middle and right, respectively). 

\subsection{Optimization and Simulation of Time Period KEPs} \label{subsec:dynamic}

As this study focuses on KEP performance differences associated with the three HLA compatibility paradigms, it relies on general and standard optimization models. These models adopt the common convention that the surgery by which recipient $r_i$ receives a transplant, say from donor $l_k$, and the surgery by which donor $l_i$ donates, say to recipient $r_s$ must occur concurrently. By consequence, these transplants form part of a cycle $C = (l_i,r_s),(l_s,r_t),\ldots, (l_t,r_k),(l_k,r_i)$, for some $i,s,t,k$ in $1,\ldots,m$ of concurrently executed transplants. We now formally define feasible cycles:

\begin{definition}
Pair $p_i$ is {\it compatible} with pair $p_j, 1 \leq i, j, \leq m$ if donor $d_i$ is compatible with recipient $r_j$, more specifically:
\begin{enumerate}
    \item $d_i$ is ABO (blood type) compatible with $r_j$  
    \item $d_i$ is HLA compatible with $r_j$ (see Section \ref{sec:casestudy} for details),
    \item $d_i$ is time compatible with $r_j$, i.e. $[a_j,d_j] \cap [a_i,d_i] \neq \emptyset$.  
\end{enumerate} \label{def:comp}
\end{definition}

\begin{definition}
A cycle $C$ of pairs $(i_1,i_2,\ldots,i_{|C|})$, where $i_1,i_2,\ldots,i_{|C|}$ in $1,\ldots,m$ is compatible if pair $i_k$ is compatible with pair $i_{k+1}$ for $k=1,\ldots,|C|-1$, pair $i_{|C|}$ is compatible with pair $i_{1}$, and $\cap_{k=1,\ldots,|C|} [a_{i_k},d_{i_k}] \neq \emptyset$.
\end{definition}

We define for any time $t, T_0 \leq t \leq T_E$, the set $\mathcal{C}_t$ of compatible cycles at time $t$, i.e., the set of compatible cycles $C = (i_1,i_2,\ldots,i_{|C|})$ such that $a_{i_k} \leq t$ and $d_{i_k} \geq t, k = 1,\ldots,|C|$. 


In practice, KEPs commonly limit the maximum cycle length. The UK KEP and some others pose a maximum length of three \citep{biro2021modelling}. Increasing the maximum cycle length beyond three may only result in marginal improvements when optimizing over a time period \citep{carvalho2023penalties}. In the remainder, we confine the analysis to a basic KEP with a maximum cycle length of three and in which the enrollment is restricted of pairs with a recipients and a single donor, while focusing on the HLA compatibility paradigms. 

For $j= 1,\dots,n$, let $t_j$ be the {\it time instants} at which cycles are formed in which recipients are matched for transplantation. Thus, $t_j$ in $\mathcal{T}$ for $j= 1,\ldots,n$ and we let $t_j < t_{j+1}, j= 1,\ldots,n-1$. The well-known cycle formulation for our problem of interest to address the first two research questions then reads \citep{abraham2007clearing}:

\begin{definition} \label{periodKEP}
For a given time period $\mathcal{T}$, set of time instants $\mathcal{I}= (t_1, t_2,\ldots,t_n)$, and set of pairs $\mathcal{P}$, the {\it time period KEP optimization problem} KEP($\mathcal{T}$) is defined as:
\begin{eqnarray}
\text{Max} \sum_{j=1}^n \sum_{C \in \mathcal{C}_{t_j}} w_{C}^{t_j} x_{C}^{t_j} & & \\
\text{Subject to} & & \\
\sum_{j=1}^n \sum_{(C \in \mathcal{C}_{t_j} | p_i \in C)} x_{C}^{t_j} & \leq & 1 \hspace{2cm} \forall i \in 1,\ldots,m \\
x_{C}^{t_j} & \in & \{0,1\} \hspace{2cm} \forall t_j \in I, C \in \mathcal{C}_j . 
\end{eqnarray}
\end{definition}

The superscript $t_j$ enforces that the time of selection is specified for each selected cycle as will be convenient in the remainder. 
The cycle weights in the objective function will be defined to reflect the sum of the HLA compatibilities between subsequent pairs in a cycle. \\ 

A main difficulty in the optimization of KEPs over time is that the arrival times $a_i$ and departure times $l_i$ of KEP($\mathcal{T}$) are only known stochastically or not all. Thus, optimality of solutions can only be ascertained by the end of a time period $\mathcal{T}$, which is obviously too late for their implementation. Any matching decisions taken without being informed about future times at which pairs arrive or leave can be non-optimal. Simulation studies on the Dutch KEP showed that consecutive time instant KEP optimization yielded only 60 percent of the post-hoc optimum when optimizing expected QALYs gained \cite{glorie2022health}  . 

Nevertheless, in reality KEPs necessarily seek to optimize the performance over time by optimizing some weighted sum of cycles at every time instant $t_j, j =1,\ldots,n$. The optimization of such {\it time instant} KEPs has also been the main problem of interest in the scientific literature on KEP optimization \citep{barkel2025operational}.

\begin{definition} \label{instantKEP}
For a given time time instant $t_j$ in $\mathcal{T}$ and set of pairs $\mathcal{P}$ matched prior to $t_j$ specified by the $x_{C}^{t_k}, k = 1,\ldots,j-1$ determined at earlier time instants, the {\it time instant KEP optimization problem} KEP($t_j$) is defined as:
\begin{eqnarray}
\text{Max} \sum_{C \in \mathcal{C}_{t_j}} w_{C}^{t_j} x_{C}^{t_j}  & & \\
\text{Subject to} & & \\
\sum_{k=1}^j \sum_{(C \in \mathcal{C}_{t_k} | p_i \in C)} x_{C}^{t_k} & \leq & 1 \hspace{2cm} \forall i \in 1,\ldots,m \\
x_{C}^{t_j} & \in & \{0,1\} \hspace{2cm} \forall C \in \mathcal{C}_{t_j}. \\
\end{eqnarray}
\end{definition}

Notice that the feasible region ensures that cycles containing a pair $p_i, i = 1,\ldots,m$ for some time instant $t_k, k = 1,\ldots,j-1$ cannot be selected at time instant $t_j$.\\

A time period KEP($\mathcal{T}$) can thus be viewed as a series of time instant KEP($t_j$)'s, $j=1,\ldots,n$ in which recipient donor pairs $p_i, i=1,\ldots,m$ which have already arrived but have not yet left nor have been matched in previous KEP($t_k$)'s, $k=1,\ldots,j-1$, can be matched. Discrete event simulation can mimic arrivals and departures of pairs and implementation of the solutions of the time instant KEP($t_j$)'s occurring over a time period $\mathcal{T}$ to evaluate the KEP performance over this period (see \citep{glorie2022health, carvalho2023penalties} for extensive examples). The next subsection elaborates how to incorporate HLA compatibility paradigms in the objective functions and constraints of time instant KEP models within this discrete event simulation approach to evaluating KEP performance over time.

\subsection{HLA Compatibility based Models to Maximize Effectiveness and Equity} \label{subsec:EquityWeighting}

A recipient may be considered insufficiently compatible with some candidate donors. This reflected in the compatibility graph by deleting arcs between pairs $(p_i,p_j)$ for which $r_i$ has DSAs against the HLA of donor $d_j$, $i,j = 1,\ldots,m$, and arcs between pairs $(p_i,p_j)$ for which the mismatch beween the HLA profile of $r_i$ and  $d_j$, $i,j = 1,\ldots,m$ exceeds a threshold, such as the 'maximum mismatched antigens' threshold adopted by OTPN \citep{optn_policy5_2014}. For antigen and allele level data, this threshold can be any number between one and the number of loci considered, e.g. six or eleven. For eplet level data this threshold can be any value between 1 and the maximum eplet mismatch load, e.g. 138.\\ 

The remainder of this subsection discusses how to consider HLA compatibility in the objective function. We adopt a weighted objective function of the total number of transplants and the HLA compatibility score. Moreover, the weight of the second component is such that the total HLA compatibility score of a solution is always less than one (transplant). This is equivalent to a hierarchical objective function that prioritizes the optimization of the number of transplants over HLA compatibility. In order to model the problem as a maximization problem, the antigen, allele, and eplet based HLA compatibility of a pair is expressed as the difference between the corresponding maximum mismatch score and the mismatch load of the pair, as also adopted in Figure \ref{fig:paragimcorrelations}.

The normalized HLA compatibility of a cycle $C = (p_{c_1},\ldots,p_{c_{|C|}})$ is then defined as follows. For any two pairs $i,j = 1,\ldots,m$, let $HLA(i,j)$ be the HLA compatibility score between donor $d_i$ and recipient $r_j$ and let $Z$ be the maximum compatibility score between a recipient and a donor (e.g. 6, 10, or 138). Then, the total normalized compatibility score of a cycle can be defined as: 
\begin{equation} \label{cycleweight}
  HLA(C) = \sum_{i=1}^{|C|} \frac{HLA(c_i,c_j)}{Z}, 
\end{equation} 
where $c_{|C|+1} = c_1$

In the remainder, KEP effectiveness is then determined by an objective function of which the first component is formed by the fraction $F(\mathcal{P})$ of pairs matched and the second component by the normalized HLA compatibility. In terms of the models of Definitions \ref{periodKEP} and \ref{instantKEP}, this boils down to letting 
\begin{equation} 
     w_{C}^{t_j} = \frac{|C| + \frac{1}{m} HLA(C)}{m}.
     \label{eq:normalizedobj}
\end{equation}
 
 The first research question can now be answered by simulations in which a time period KEP($\mathcal{T}$) is solved by solving the time instant KEPs that arise over the course of time period $\mathcal{T}$ with an objective function of the form (\ref{eq:normalizedobj}). \\

Maximizing effectiveness according to (\ref{eq:normalizedobj}) and (\ref{cycleweight}) may be more favorable for some patient subpopulations than for others in optimal time instant KEP solutions as may subsequently reflect in the performance of the time period KEPs. More specifically, it has been conjectured that basing HLA compatibility on allele level data rather than on antigen level data reduces HLA compatibility scores for ethnic minorities and thus induces or enlarges differences in $F(s)$, $L(s)$, and $HLA(s)$ of (minority) subpopulations $s$ with the overall $F(P)$, $L(P)$, and $HLA(P)$.

In pursuit of answering the second research question, equity is therefore defined in terms of differences in the objective function values (\ref{eq:normalizedobj}) obtained for subpopulations $s$ with the values obtained for the general population:
\begin{equation} 
     \frac{|s|}{m} [ (F(s) - F(P))  + \frac{1}{m} ( HLA(s) - HLA(P) ) ]
\label{eq:equity}
\end{equation}



To answer the third research question, we adopt an equity weighting approach to jointly optimize effectiveness and equity \cite{cookson2020distributional}. The approach minimizes the objective function of the time instant $KEP(t_j)$s that arise in KEP($\mathcal{T}$) using time-independent subpopulation weights $v_s, s \in S$ as follows. For $s \in S$, let $i=1,\ldots,m$, let $s(r_i)$ be the subpopulation to which recipient $r_i$ belongs. A weighted objective function for a time instant KEP is now (re)defined as:\\   
\begin{equation} \label{weightedobjfunction}
     w_{C}^{t_j} = \sum_{i=1}^{|C|} v(s(r_i)) \times [ 1 +  \frac{HLA(c_i,c_j)}{m \times Z}]. 
\end{equation}

The subsequent question is to determine weights that eliminate inequities in (\ref{eq:equity}). The differences achieved for a time period KEP are not a direct optimization result but determined by optimally solving a series of time instant KEPs.  In pursuit of equity, one may seek weights that jointly minimize: 
\begin{equation} 
     \sum_{s \in T}  \| \frac{|s|}{m} [ (F(s) - F(P))  + \frac{1}{m} ( HLA(s) - HLA(P) ) \| ]
     \label{eq:equityweighting}
\end{equation}
where $T$ is the set of subpopulations and the sum is over the absolute differences. In the case study, weights that minimize (\ref{eq:equityweighting}) are determined following a Rawlsian approach, recursively increasing the weight of a subpopulation $s$ with the lowest objective function value score $(F(s) + \frac{1}{m} HLA(s))$ until this causes a subpopulation for which the weight has been increased to have a larger positive difference with $(F(P) + \frac{1}{m} HLA(P))$ than any remaining negative differences for other subpopulations \citep{rawls2017theory}. In view of its heuristic and simulation-based nature, the set of weights with lowest value for (\ref{eq:equityweighting}) encountered in the process is reported as most equitable.

\section{Case Study} \label{sec:casestudy}

\subsection{Case Study Data Collection and Calibration} \label{subsec:data}

The case study is based on a cohort of recipients and donors extracted from the OPTN star file under a data sharing agreement and an accompanying data set on donor specific antigens \citep{optndataset2024}. Data were extracted from the 1332 recipients and 1401 donors for which allele-level HLA data for each of the HLA A, B and DR loci were available (as needed for imputation, see \ref{subsubsec:compatibilityGraphs}. From these two sets, a maximum number of 990 of mutually incompatible recipient-donor pairs was formed by solving a maximum cardinality bipartite matching problem. Thus, the set of pairs forming the KEP population of the case study is constructed and not representative of an actual KEP. \\

Each simulation run covers a period of ten years in which a time instant KEP is solved every three months. The arrival times of the pairs are sampled from an exponential distribution which in expectation generates 990 arrivals in the ten years time interval. The departure time of each pair is sampled from an exponential distribution of which the departure rate is 0.29 times the arrival rate (as reported for the UNOS/OPTN waitlist \cite{us2024health}). 

The discrete event simulation is programmed in Python version 3.13.5 and the integer programming formulations of the time instant KEPs are solved using Gurobi version 12.0.3. The code is on Github 
\citep{PeraltaGithub2025}. Averages and standard deviations over 100 simulations are reported.

\subsubsection{HLA data and compatibility Graphs} \label{subsubsec:compatibilityGraphs}

The available HLA typings of recipients and donors regard the HLA-A, B, C, DR, and DQ loci. No data were available for the DP locus. With two antigens (alleles) per loci, this results in a maximum of 10 matches. For the recipients and donors for whom there were no allele level data for all loci considered, allele level data were imputed based on the allele level HLA A,B, and DR data, any additional allele or antigen level data available, and ethnicity using Haplostats \cite{gragert2013six}. Haplostats requires HLA A,B, and DR data to be specified at the allele level. All selected donors had antigen or allele level data for the additional loci HLA-C and DQ. For some of the selected recipients, the HLA-C or DQ were missing or incomplete.

DSA data were reported for most recipients, either at the antigen or at the allele level. DSA based incompatibility was assessed using the original (non-imputed) HLA typings taking a conservative approach. A pair was classified as incompatible if (i) the DSA matched the donor at the exact allele level, (ii) an antigen-level DSA matched the donor’s antigen, or (iii) an allele-level DSA mapped to the donor’s reported antigen family, or (iv) an antigen-level DSA matched the antigen family of a donor allele reported at allele level. The following example illustrates this conservative approach using antigen/allele B07/B*07:xx):

\begin{enumerate}
    \item Donor \textbf{B07}; patient DSA \textbf{B*07:01}: incompatible.
    \item Donor \textbf{B07}; patient DSA \textbf{B07}: incompatible.
    \item Donor \textbf{B*07:01}; patient DSA \textbf{B07}: incompatible.
    \item Donor \textbf{B*07:01}; patient DSA \textbf{B*07:01}: incompatible.
    \item Donor \textbf{B*07:01}; patient DSA \textbf{B*07:02}: compatible.
\end{enumerate}

The eplet match score was computed following the definitions of the HLA Eplet Registry \citep{HLAEpletRegistry2025}. The implementation requires high resolution HLA typings for donors and recipients at the HLA-A, -B, -C, -DR, and -DQ loci and reports separate eplet match scores for DR, DQ, and for the three class 1 loci HLA-A, -B, and -C together. When the original HLA data were provided at the allele level, we mapped alleles to their corresponding antigens following Eurotransplant HLA tables to calculate antigen level HLA scores \citep{ETRLHLAtable2023}.   

Based on the resulting mismatch scores for each of the three compatibility paradigms, a compatibility graph was subsequently constructed following Definition \ref{def:comp}. In further specification of this definition, a recipient and donor pair is considered HLA compatible if: 

\begin{enumerate}
    \item no donor-specific antibody (DSA) are present against any of the donor’s HLA antigens or alleles,
    \item the HLA match score satisfied the minimum compatibility thresholds of 3, 2 and 82 for antigens, alleles, and eplets respectively, based on the ten loci HLA profiles. 
 \end{enumerate}
These thresholds ensure that recipients are not (unrealistically) matched with any donor for which they don´t have DSAs. Moreover, the threshold values yield comparable numbers of arcs in the resulting compatibility graphs of 251.863, 216.604, and 234.601 respectively. 

For the simulation studies in which only the HLA B, DR, DQ loci are considered, the threshold values are 2 and 1 for antigen and allele level respectively, resulting in 242.708 and 274.942 respectively. Eplet load cannot be calculated for the combination of HLA-B, DR, and DQ. For DR and DQ loci simulation studies, these threshold values are 2, 1, and 45 for antigen, allele and eplet load and result in 203.551, 210.153, and 203.504 arcs respectively. Eplet load is then limited to the (Class II) DR and DQ loci.

\subsection{Simulation Results} \label{subsec:results}

Table \ref{tab:orig_weights_10loci} below presents the results for the ten year KEP when optimizing the time instant KEPs based on (\ref{eq:normalizedobj}) and using antigen-, allele-, and eplet-level data respectively. In addition to results for the total population of patients participating in the KEP, it presents subpopulation results for each of the ethnic subpopulations distinguished by UNOS. For each (sub)population $s$ and for the entire population $P$, the results include the 100 simulation averages of the fraction $F(s)$ of patients matched, the fraction $L(s)$ of patients that have left without being matched, and the remaining fraction of patients who are still in the KEP by $T_E$. Moreover, it presents the HLA match scores $HLA(s)$, as well as 95 percent confidence intervals for each of the compatibility paradigms. The average waiting times $W(s)$, expressed in number of KEPs participated before being matched, is reported as well.

\begin{table}[htbp]
\centering
\caption{Unweighted maximizing of the number of transplants (HLA scores based on 10 loci)}
\label{tab:orig_weights_10loci}
\captionsetup[subtable]{labelformat=empty,justification=centering}
\footnotesize
\setlength{\tabcolsep}{2.6pt}
\renewcommand{\arraystretch}{1.20}

\begin{subtable}{\textwidth}
\centering
\caption{\textbf{Antigen level based optimization and compatibility}}
\label{tab:antigen_level_10loci}
\begin{adjustbox}{max width=\textwidth}
\begin{tabular}{
  L{3.4cm} C{1.6cm} C{2.2cm} C{2.7cm} C{2.7cm} C{3.4cm} C{1.8cm} C{2.35cm} C{2.60cm}}
\toprule
\multicolumn{1}{c}{\textbf{Ethnicity}} & \multicolumn{1}{c}{\textbf{Arrivals}} &
\multicolumn{1}{c}{\(\bm{F(s)}\)} & \multicolumn{3}{c}{\textbf{HLA(s)}} &
\multicolumn{1}{c}{\(\bm{W(s)}\)} & \multicolumn{1}{c}{\(\bm{L(s)}\)} &
\multicolumn{1}{c}{\(\bm{1 - F(s) - L^{\ast}(s)}\)} \\
\cmidrule(lr){4-6}
\multicolumn{1}{c}{} & \multicolumn{1}{c}{} &
\multicolumn{1}{c}{\footnotesize(Matched)} &
\multicolumn{1}{c}{\textbf{Antigen}} &
\multicolumn{1}{c}{\textbf{Allele}} &
\multicolumn{1}{c}{\textbf{Eplets}} &
\multicolumn{1}{c}{\footnotesize(Waiting Time)} &
\multicolumn{1}{c}{\footnotesize(Left Unmatched)} &
\multicolumn{1}{c}{\footnotesize(Still in KEP)} \\
\midrule
Caucasian & 592.62 & \makecell{\underline{0.699}\\ \underline{\ci{0.695}{0.703}}} &
\makecell{\underline{4.705}\\ \underline{\ci{4.690}{4.720}}} &
\makecell{\underline{3.399}\\ \underline{\ci{3.381}{3.417}}} &
\makecell{\underline{93.181}\\ \underline{\ci{93.047}{93.314}}} &
2.977 & \makecell{\underline{0.269}\\ \underline{\ci{0.265}{0.273}}} & 0.032 \\
Afroamerican & 170.18 & \makecell{\underline{0.636}\\ \underline{\ci{0.629}{0.643}}} &
\makecell{\underline{4.150}\\ \underline{\ci{4.128}{4.173}}} &
\makecell{\underline{2.537}\\ \underline{\ci{2.511}{2.563}}} &
\makecell{\underline{90.574}\\ \underline{\ci{90.362}{90.787}}} &
3.109 & \makecell{\underline{0.325}\\ \underline{\ci{0.318}{0.331}}} & 0.039 \\
Latin & 147.40 & \makecell{\underline{0.622}\\ \underline{\ci{0.615}{0.629}}} &
\makecell{\underline{4.293}\\ \underline{\ci{4.267}{4.320}}} &
\makecell{\underline{2.656}\\ \underline{\ci{2.628}{2.685}}} &
\makecell{\underline{91.005}\\ \underline{\ci{90.754}{91.257}}} &
3.533 & \makecell{\underline{0.334}\\ \underline{\ci{0.327}{0.341}}} & 0.044 \\
Asian & 55.34 & \makecell{\underline{0.627}\\ \underline{\ci{0.616}{0.639}}} &
\makecell{\underline{3.934}\\ \underline{\ci{3.901}{3.966}}} &
\makecell{\underline{2.059}\\ \underline{\ci{2.020}{2.099}}} &
\makecell{\underline{84.786}\\ \underline{\ci{84.367}{85.205}}} &
3.737 & \makecell{\underline{0.330}\\ \underline{\ci{0.317}{0.342}}} & 0.043 \\
American ind.\ / Alaska nat. & 7.95 &
\makecell{\underline{0.763}\\ \underline{\ci{0.736}{0.789}}} &
\makecell{\underline{4.219}\\ \underline{\ci{4.119}{4.320}}} &
\makecell{\underline{2.923}\\ \underline{\ci{2.801}{3.044}}} &
\makecell{\underline{84.346}\\ \underline{\ci{83.381}{85.310}}} &
1.979 & \makecell{\underline{0.213}\\ \underline{\ci{0.186}{0.241}}} & 0.024 \\
Native Hawaiian\ / Pacific islander & 5.94 &
\makecell{\underline{0.461}\\ \underline{\ci{0.428}{0.494}}} &
\makecell{\underline{4.101}\\ \underline{\ci{3.964}{4.239}}} &
\makecell{\underline{2.266}\\ \underline{\ci{2.105}{2.427}}} &
\makecell{\underline{92.870}\\ \underline{\ci{91.145}{94.595}}} &
5.668 & \makecell{\underline{0.468}\\ \underline{\ci{0.433}{0.503}}} & 0.071 \\
\midrule
\textbf{Entire Population} & \textbf{979.43} &
\makecell{\textbf{0.672}\\ \footnotesize\textbf{[0.668;0.675]}} &
\makecell{\textbf{4.508}\\ \footnotesize\textbf{[4.497; 4.520]}} &
\makecell{\textbf{3.074}\\ \footnotesize\textbf{[3.061; 3.086]}} &
\makecell{\textbf{91.921}\\ \footnotesize\textbf{[91.822; 92.019]}} &
\textbf{3.120} & \makecell{\textbf{0.293}\\ \footnotesize\textbf{[0.289; 0.296]}} & \textbf{0.035} \\
\bottomrule
\end{tabular}
\end{adjustbox}
\end{subtable}

\vspace{0.8em}

\begin{subtable}{\textwidth}
\centering
\caption{\textbf{Allele level based optimization and compatibility}}
\label{tab:allele_level_10loci}
\begin{adjustbox}{max width=\textwidth}
\begin{tabular}{
  L{3.4cm} C{1.6cm} C{2.2cm} C{2.7cm} C{2.7cm} C{3.4cm} C{1.8cm} C{2.35cm} C{2.60cm}}
\toprule
\multicolumn{1}{c}{\textbf{Ethnicity}} & \multicolumn{1}{c}{\textbf{Arrivals}} &
\multicolumn{1}{c}{\(\bm{F(s)}\)} & \multicolumn{3}{c}{\textbf{HLA(s)}} &
\multicolumn{1}{c}{\(\bm{W(s)}\)} & \multicolumn{1}{c}{\(\bm{L(s)}\)} &
\multicolumn{1}{c}{\(\bm{1 - F(s) - L^{\ast}(s)}\)} \\
\cmidrule(lr){4-6}
\multicolumn{1}{c}{} & \multicolumn{1}{c}{} &
\multicolumn{1}{c}{\footnotesize(Matched)} &
\multicolumn{1}{c}{\textbf{Antigen}} &
\multicolumn{1}{c}{\textbf{Allele}} &
\multicolumn{1}{c}{\textbf{Eplets}} &
\multicolumn{1}{c}{\footnotesize(Waiting Time)} &
\multicolumn{1}{c}{\footnotesize(Left Unmatched)} &
\multicolumn{1}{c}{\footnotesize(Still in KEP)} \\
\midrule
Caucasian & 595.14 & \makecell{\underline{0.704}\\ \underline{\ci{0.700}{0.708}}} &
\makecell{\underline{4.664}\\ \underline{\ci{4.650}{4.678}}} &
\makecell{\underline{3.648}\\ \underline{\ci{3.633}{3.662}}} &
\makecell{\underline{93.038}\\ \underline{\ci{92.912}{93.163}}} &
3.33 & \makecell{\underline{0.262}\\ \underline{\ci{0.257}{0.266}}} & 0.035 \\
Afroamerican & 170.9 & \makecell{\underline{0.635}\\ \underline{\ci{0.629}{0.642}}} &
\makecell{\underline{4.022}\\ \underline{\ci{3.995}{4.049}}} &
\makecell{\underline{2.954}\\ \underline{\ci{2.937}{2.970}}} &
\makecell{\underline{90.074}\\ \underline{\ci{89.830}{90.317}}} &
3.71 & \makecell{\underline{0.318}\\ \underline{\ci{0.312}{0.325}}} & 0.046 \\
Latin & 147.93 & \makecell{\underline{0.614}\\ \underline{\ci{0.606}{0.621}}} &
\makecell{\underline{4.254}\\ \underline{\ci{4.226}{4.281}}} &
\makecell{\underline{3.065}\\ \underline{\ci{3.043}{3.087}}} &
\makecell{\underline{91.401}\\ \underline{\ci{91.147}{91.656}}} &
4.018 & \makecell{\underline{0.337}\\ \underline{\ci{0.330}{0.344}}} & 0.049 \\
Asian & 55.63 & \makecell{\underline{0.535}\\ \underline{\ci{0.525}{0.546}}} &
\makecell{\underline{3.900}\\ \underline{\ci{3.859}{3.942}}} &
\makecell{\underline{2.633}\\ \underline{\ci{2.599}{2.667}}} &
\makecell{\underline{85.068}\\ \underline{\ci{84.597}{85.539}}} &
5.182 & \makecell{\underline{0.401}\\ \underline{\ci{0.390}{0.412}}} & 0.064 \\
American ind.\ / Alaska nat. & 7.90 &
\makecell{\underline{0.824}\\ \underline{\ci{0.798}{0.850}}} &
\makecell{\underline{4.088}\\ \underline{\ci{3.975}{4.202}}} &
\makecell{\underline{3.164}\\ \underline{\ci{3.081}{3.247}}} &
\makecell{\underline{83.989}\\ \underline{\ci{82.859}{85.120}}} &
1.99 & \makecell{\underline{0.159}\\ \underline{\ci{0.133}{0.185}}} & 0.017 \\
Native Hawaiian\ / Pacific islander & 5.96 &
\makecell{\underline{0.318}\\ \underline{\ci{0.287}{0.348}}} &
\makecell{\underline{3.952}\\ \underline{\ci{3.752}{4.152}}} &
\makecell{\underline{2.957}\\ \underline{\ci{2.803}{3.111}}} &
\makecell{\underline{94.330}\\ \underline{\ci{92.806}{95.853}}} &
6.755 & \makecell{\underline{0.566}\\ \underline{\ci{0.529}{0.602}}} & 0.117 \\
\midrule
\textbf{Entire Population} & \textbf{983.46} &
\makecell{\textbf{0.667}\\ \footnotesize\textbf{[0.664;0.671]}} &
\makecell{\textbf{4.459}\\ \footnotesize\textbf{[4.448; 4.470]}} &
\makecell{\textbf{3.399}\\ \footnotesize\textbf{[3.388; 3.411]}} &
\makecell{\textbf{91.874}\\ \footnotesize\textbf{[91.774; 91.974]}} &
\textbf{3.574} & \makecell{\textbf{0.292}\\ \footnotesize\textbf{[0.288; 0.295]}} & \textbf{0.041} \\
\bottomrule
\end{tabular}
\end{adjustbox}
\end{subtable}

\vspace{0.8em}

\begin{subtable}{\textwidth}
\centering
\caption{\textbf{Eplet level based optimization and compatibility}}
\label{tab:eplet_level_10loci}
\begin{adjustbox}{max width=\textwidth}
\begin{tabular}{
  L{3.4cm} C{1.6cm} C{2.2cm} C{2.7cm} C{2.7cm} C{3.4cm} C{1.8cm} C{2.35cm} C{2.60cm}}
\toprule
\multicolumn{1}{c}{\textbf{Ethnicity}} & \multicolumn{1}{c}{\textbf{Arrivals}} &
\multicolumn{1}{c}{\(\bm{F(s)}\)} & \multicolumn{3}{c}{\textbf{HLA(s)}} &
\multicolumn{1}{c}{\(\bm{W(s)}\)} & \multicolumn{1}{c}{\(\bm{L(s)}\)} &
\multicolumn{1}{c}{\(\bm{1 - F(s) - L^{\ast}(s)}\)} \\
\cmidrule(lr){4-6}
\multicolumn{1}{c}{} & \multicolumn{1}{c}{} &
\multicolumn{1}{c}{\footnotesize(Matched)} &
\multicolumn{1}{c}{\textbf{Antigen}} &
\multicolumn{1}{c}{\textbf{Allele}} &
\multicolumn{1}{c}{\textbf{Eplets}} &
\multicolumn{1}{c}{\footnotesize(Waiting Time)} &
\multicolumn{1}{c}{\footnotesize(Left Unmatched)} &
\multicolumn{1}{c}{\footnotesize(Still in KEP)} \\
\midrule
Caucasian & 590.98 & \makecell{\underline{0.684}\\ \underline{\ci{0.680}{0.688}}} &
\makecell{\underline{4.610}\\ \underline{\ci{4.596}{4.624}}} &
\makecell{\underline{3.345}\\ \underline{\ci{3.330}{3.360}}} &
\makecell{\underline{96.742}\\ \underline{\ci{96.638}{96.846}}} &
3.43 & \makecell{\underline{0.268}\\ \underline{\ci{0.264}{0.272}}} & 0.048 \\
Afroamerican & 169.76 & \makecell{\underline{0.632}\\ \underline{\ci{0.626}{0.638}}} &
\makecell{\underline{3.902}\\ \underline{\ci{3.874}{3.930}}} &
\makecell{\underline{2.380}\\ \underline{\ci{2.353}{2.407}}} &
\makecell{\underline{93.909}\\ \underline{\ci{93.736}{94.081}}} &
3.656 & \makecell{\underline{0.304}\\ \underline{\ci{0.297}{0.311}}} & 0.064 \\
Latin & 147.12 & \makecell{\underline{0.591}\\ \underline{\ci{0.585}{0.597}}} &
\makecell{\underline{4.148}\\ \underline{\ci{4.114}{4.181}}} &
\makecell{\underline{2.590}\\ \underline{\ci{2.562}{2.618}}} &
\makecell{\underline{94.718}\\ \underline{\ci{94.496}{94.939}}} &
4.328 & \makecell{\underline{0.336}\\ \underline{\ci{0.329}{0.342}}} & 0.073 \\
Asian & 55.24 & \makecell{\underline{0.541}\\ \underline{\ci{0.529}{0.552}}} &
\makecell{\underline{3.923}\\ \underline{\ci{3.872}{3.974}}} &
\makecell{\underline{1.964}\\ \underline{\ci{1.918}{2.010}}} &
\makecell{\underline{92.642}\\ \underline{\ci{92.296}{92.988}}} &
4.895 & \makecell{\underline{0.372}\\ \underline{\ci{0.361}{0.384}}} & 0.087 \\
American ind.\ / Alaska nat. & 7.91 &
\makecell{\underline{0.697}\\ \underline{\ci{0.671}{0.724}}} &
\makecell{\underline{4.530}\\ \underline{\ci{4.414}{4.645}}} &
\makecell{\underline{3.118}\\ \underline{\ci{2.991}{3.246}}} &
\makecell{\underline{92.997}\\ \underline{\ci{92.233}{93.761}}} &
2.801 & \makecell{\underline{0.263}\\ \underline{\ci{0.236}{0.291}}} & 0.039 \\
Native Hawaiian\ / Pacific islander & 5.91 &
\makecell{\underline{0.497}\\ \underline{\ci{0.471}{0.524}}} &
\makecell{\underline{3.536}\\ \underline{\ci{3.381}{3.690}}} &
\makecell{\underline{2.080}\\ \underline{\ci{1.938}{2.222}}} &
\makecell{\underline{94.874}\\ \underline{\ci{93.842}{95.906}}} &
4.669 & \makecell{\underline{0.411}\\ \underline{\ci{0.378}{0.445}}} & 0.091 \\
\midrule
\textbf{Entire Population} & \textbf{976.92} &
\makecell{\textbf{0.652}\\ \footnotesize\textbf{[0.649;0.655]}} &
\makecell{\textbf{4.390}\\ \footnotesize\textbf{[4.379; 4.401]}} &
\makecell{\textbf{3.007}\\ \footnotesize\textbf{[2.995; 3.018]}} &
\makecell{\textbf{95.754}\\ \footnotesize\textbf{[95.671; 95.836]}} &
\textbf{3.665} & \makecell{\textbf{0.291}\\ \footnotesize\textbf{[0.288; 0.295]}} & \textbf{0.057} \\
\bottomrule
\end{tabular}
\end{adjustbox}
\end{subtable}

\end{table}

For each metric, the performance for a subpopulation is considered to differ significantly from the overall performance if the subpopulation average is outside of the 95 percent confidence interval for the overall population, based on the observed mean and variance and assuming a normal distribution. All significant results are underlined and one easily verifies that almost all subpopulation results differ significantly from the overall results. While the transplant probabilities for the entire population are close to two-third for each of the three paradigms, these probabilities range from close to 0.8 to below 0.4 for some of the subpopulations. 

Caucasians consistently have significantly higher transplant probabilities and HLA match scores than most other subpopulations while the opposite holds for Asians. The subpopulations of American indians / Alaskan natives and of Native Hawaiian / Pacific islanders are so small that they are not represented in the majority of time instant KEPs and the variances in performance for these subpopulations are large. Given their small sizes, it appears more appropriate to consider effectiveness and efficiency for these minorities using other evaluation methods and the remainder of the analysis focuses on the other four subpopulations.  

In general, the probabilities of receiving a transplant appear to be positively correlated with HLA scores, as likely caused by the second component in the objective function. When comparing results for a single subpopulation across the HLA paradigms, the differences in transplant probabilities are very small and mostly non-significant. For the Asian subpopulation, however, the probabilities obtained when using allele and eplet level data yield significantly and substantially lower transplant probabilities. This is accompanied by longer waiting times to transplant and higher probabilities of leaving the KEP without a transplant when using allele and eplet paradigms compared to antigen based HLA compatibility. 

For the entire population, the average antigen level HLA match obtained with the antigen paradigm is around 4.5, almost half of the maximum match score of 10. The same paradigm yields much lower results at the allele level, around 3.1 and an average eplet score of 91.9 (for which the maximum is 138). These results change only for the allele level score when using allele level data for compatibility and optimization and the allele level score increases to 3.4. Eplet based optimization increase the eplet score significantly to 95.7 and lowers the average antigen and allele level scores to 4.4 and 3.0 respectively. 

For the constructed KEP in the case study, the results show that changing the HLA compatibility paradigms in KEP feasibility and optimization severely impacts effectiveness and equity. Moreover, accuracy gains made possible by the adoption of the allele and eplet level compatibility paradigms  negatively impacted equity. The supplementary file presents results for different combinations of feasibility and optimization paradigms, which largely confirm the aforementioned results and show that the results are not caused by differences in compatibilty graphs (feasible solution space). 

The supplementary file also presents result tables for HLA compatibilty scores based on the HLA B, DR, and DQ loci and on the DR and DQ loci only. Eplet mismatch cannot be calculated for the former.

\subsection{Equity Improvement} \label{subsec:EquityImprovement}

The third research question was answered by adjusting the subpopulation weights in (\ref{weightedobjfunction}) with small increments of 0.001.  Table \ref{tab:equity_weighting_10loci} presents the weights and metrics of KEP performance when considering 10 loci (4 and 6 loci results are presented in the supplementary materials).

\begin{table}[htbp]
\centering
\caption{Maximizing the equity weighted number of transplants (HLA scores based on 10 loci)}
\label{tab:equity_weighting_10loci}
\captionsetup[subtable]{labelformat=empty,justification=centering}
\footnotesize
\setlength{\tabcolsep}{2.4pt}
\renewcommand{\arraystretch}{1.20}

\begin{subtable}{\textwidth}
\centering
\caption{\textbf{Antigen level based optimization and compatibility}}
\label{tab:lr_lr_mod_sub}
\begin{adjustbox}{max width=\textwidth}
\begin{tabular}{
  L{3.4cm}  
  C{2.0cm}  
  C{1.6cm}  
  C{2.2cm}  
  C{2.7cm}  
  C{2.7cm}  
  C{3.4cm}  
  C{1.8cm}  
  C{2.35cm} 
  C{2.60cm} 
}
\toprule
\multicolumn{1}{c}{\textbf{Ethnicity}} &
\multicolumn{1}{c}{\textbf{Equity weight}} &
\multicolumn{1}{c}{\textbf{Arrivals}} &
\multicolumn{1}{c}{\(\bm{F(s)}\)} &
\multicolumn{3}{c}{\textbf{HLA(s)}} &
\multicolumn{1}{c}{\(\bm{W(s)}\)} &
\multicolumn{1}{c}{\(\bm{L(s)}\)} &
\multicolumn{1}{c}{\(\bm{1 - F(s) - L^{\ast}(s)}\)} \\
\cmidrule(lr){5-7}
\multicolumn{1}{c}{} & \multicolumn{1}{c}{} & \multicolumn{1}{c}{} &
\multicolumn{1}{c}{\footnotesize(Matched)} &
\multicolumn{1}{c}{\textbf{Antigen}} &
\multicolumn{1}{c}{\textbf{Allele}} &
\multicolumn{1}{c}{\textbf{Eplets}} &
\multicolumn{1}{c}{\footnotesize(Waiting Time)} &
\multicolumn{1}{c}{\footnotesize(Left Unmatched)} &
\multicolumn{1}{c}{\footnotesize(Still in KEP)} \\
\midrule
Caucasian & 1.000 & 591.56 &
\makecell{0.673\\ \ci{0.669}{0.677}} &
\makecell{\underline{4.683}\\ \underline{\ci{4.668}{4.697}}} &
\makecell{\underline{3.394}\\ \underline{\ci{3.379}{3.410}}} &
\makecell{\underline{93.036}\\ \underline{\ci{92.907}{93.165}}} &
3.095 &
\makecell{0.290\\ \ci{0.287}{0.294}} &
0.036 \\
Afroamerican & 1.003 & 170.05 &
\makecell{\underline{0.679}\\ \underline{\ci{0.672}{0.686}}} &
\makecell{\underline{4.174}\\ \underline{\ci{4.148}{4.199}}} &
\makecell{\underline{2.534}\\ \underline{\ci{2.509}{2.558}}} &
\makecell{\underline{90.787}\\ \underline{\ci{90.567}{91.006}}} &
3.128 &
\makecell{\underline{0.286}\\ \underline{\ci{0.279}{0.293}}} &
0.035 \\
Latin & 1.003 & 147.22 &
\makecell{\underline{0.665}\\ \underline{\ci{0.658}{0.672}}} &
\makecell{\underline{4.324}\\ \underline{\ci{4.298}{4.350}}} &
\makecell{\underline{2.690}\\ \underline{\ci{2.663}{2.717}}} &
\makecell{\underline{91.165}\\ \underline{\ci{90.896}{91.433}}} &
3.264 &
\makecell{\underline{0.299}\\ \underline{\ci{0.291}{0.307}}} &
0.036 \\
Asian & 1.003 & 55.31 &
\makecell{\underline{0.667}\\ \underline{\ci{0.654}{0.680}}} &
\makecell{\underline{3.935}\\ \underline{\ci{3.900}{3.969}}} &
\makecell{\underline{2.033}\\ \underline{\ci{2.001}{2.065}}} &
\makecell{\underline{84.352}\\ \underline{\ci{83.947}{84.757}}} &
3.386 &
\makecell{\underline{0.299}\\ \underline{\ci{0.286}{0.312}}} &
0.033 \\
American ind.\ / Alaska nat. & 1.000 & 7.94 &
\makecell{\underline{0.740}\\ \underline{\ci{0.714}{0.767}}} &
\makecell{\underline{4.219}\\ \underline{\ci{4.121}{4.317}}} &
\makecell{\underline{2.873}\\ \underline{\ci{2.765}{2.981}}} &
\makecell{\underline{83.872}\\ \underline{\ci{82.857}{84.887}}} &
2.102 &
\makecell{\underline{0.234}\\ \underline{\ci{0.207}{0.261}}} &
0.025 \\
Native Hawaiian\ / Pacific islander & 1.000 & 5.93 &
\makecell{\underline{0.450}\\ \underline{\ci{0.415}{0.486}}} &
\makecell{\underline{4.142}\\ \underline{\ci{3.998}{4.286}}} &
\makecell{\underline{2.322}\\ \underline{\ci{2.163}{2.481}}} &
\makecell{\underline{93.955}\\ \underline{\ci{92.581}{95.329}}} &
5.63 &
\makecell{\underline{0.474}\\ \underline{\ci{0.437}{0.512}}} &
0.075 \\
\midrule
\textbf{Entire Population} & \textemdash{} & \textbf{978.01} &
\makecell{\textbf{0.672}\\ \footnotesize\textbf{[0.669;0.675]}} &
\makecell{\textbf{4.491}\\ \footnotesize\textbf{[4.481; 4.502]}} &
\makecell{\textbf{3.053}\\ \footnotesize\textbf{[3.041; 3.065]}} &
\makecell{\textbf{91.796}\\ \footnotesize\textbf{[91.698; 91.894]}} &
\textbf{3.146} &
\makecell{\textbf{0.292}\\ \footnotesize\textbf{[0.289; 0.296]}} &
\textbf{0.036} \\
\bottomrule
\end{tabular}
\end{adjustbox}
\end{subtable}

\vspace{0.8em}

\begin{subtable}{\textwidth}
\centering
\caption{\textbf{Allele level based optimization and compatibility}}
\label{tab:hr_hr_mod_sub}
\begin{adjustbox}{max width=\textwidth}
\begin{tabular}{
  L{3.4cm} C{2.0cm} C{1.6cm} C{2.2cm} C{2.7cm} C{2.7cm} C{3.4cm} C{1.8cm} C{2.35cm} C{2.60cm}}
\toprule
\multicolumn{1}{c}{\textbf{Ethnicity}} &
\multicolumn{1}{c}{\textbf{Equity weight}} &
\multicolumn{1}{c}{\textbf{Arrivals}} &
\multicolumn{1}{c}{\(\bm{F(s)}\)} &
\multicolumn{3}{c}{\textbf{HLA(s)}} &
\multicolumn{1}{c}{\(\bm{W(s)}\)} &
\multicolumn{1}{c}{\(\bm{L(s)}\)} &
\multicolumn{1}{c}{\(\bm{1 - F(s) - L^{\ast}(s)}\)} \\
\cmidrule(lr){5-7}
\multicolumn{1}{c}{} & \multicolumn{1}{c}{} & \multicolumn{1}{c}{} &
\multicolumn{1}{c}{\footnotesize(Matched)} &
\multicolumn{1}{c}{\textbf{Antigen}} &
\multicolumn{1}{c}{\textbf{Allele}} &
\multicolumn{1}{c}{\textbf{Eplets}} &
\multicolumn{1}{c}{\footnotesize(Waiting Time)} &
\multicolumn{1}{c}{\footnotesize(Left Unmatched)} &
\multicolumn{1}{c}{\footnotesize(Still in KEP)} \\
\midrule
Caucasian & 1.000 & 590.13 &
\makecell{0.662\\ \ci{0.658}{0.666}} &
\makecell{\underline{4.627}\\ \underline{\ci{4.612}{4.642}}} &
\makecell{\underline{3.604}\\ \underline{\ci{3.591}{3.618}}} &
\makecell{\underline{92.757}\\ \underline{\ci{92.627}{92.888}}} &
3.425 &
\makecell{0.295\\ \ci{0.291}{0.299}} &
0.043 \\
Afroamerican & 1.003 & 169.73 &
\makecell{0.662\\ \ci{0.655}{0.670}} &
\makecell{\underline{3.997}\\ \underline{\ci{3.971}{4.023}}} &
\makecell{\underline{2.933}\\ \underline{\ci{2.915}{2.951}}} &
\makecell{\underline{89.838}\\ \underline{\ci{89.609}{90.067}}} &
3.534 &
\makecell{0.301\\ \ci{0.294}{0.308}} &
 0.037 \\
Latin & 1.004 & 146.94 &
\makecell{\underline{0.669}\\ \underline{\ci{0.662}{0.676}}} &
\makecell{\underline{4.227}\\ \underline{\ci{4.198}{4.255}}} &
\makecell{\underline{3.038}\\ \underline{\ci{3.015}{3.062}}} &
\makecell{\underline{90.946}\\ \underline{\ci{90.688}{91.203}}} &
3.585 &
\makecell{0.294\\ \ci{0.287}{0.301}} &
0.037 \\
Asian & 1.006 & 55.25 &
\makecell{0.663\\ \ci{0.651}{0.674}} &
\makecell{\underline{3.813}\\ \underline{\ci{3.775}{3.852}}} &
\makecell{\underline{2.567}\\ \underline{\ci{2.541}{2.594}}} &
\makecell{\underline{84.528}\\ \underline{\ci{84.104}{84.953}}} &
4.133 &
\makecell{\underline{0.303}\\ \underline{\ci{0.291}{0.314}}} &
0.035 \\
American ind.\ / Alaska nat. & 1.000 & 7.88 &
\makecell{\underline{0.783}\\ \underline{\ci{0.761}{0.806}}} &
\makecell{\underline{3.995}\\ \underline{\ci{3.856}{4.133}}} &
\makecell{\underline{3.132}\\ \underline{\ci{3.029}{3.234}}} &
\makecell{\underline{82.667}\\ \underline{\ci{81.517}{83.816}}} &
2.202 &
\makecell{\underline{0.192}\\ \underline{\ci{0.169}{0.214}}} &
0.025 \\
Native Hawaiian\ / Pacific islander & 1.000 & 5.89 &
\makecell{\underline{0.267}\\ \underline{\ci{0.233}{0.302}}} &
\makecell{\underline{4.004}\\ \underline{\ci{3.800}{4.208}}} &
\makecell{\underline{2.920}\\ \underline{\ci{2.754}{3.087}}} &
\makecell{\underline{96.567}\\ \underline{\ci{94.755}{98.379}}} &
7.334 &
\makecell{\underline{0.595}\\ \underline{\ci{0.557}{0.634}}} &
0.137 \\
\midrule
\textbf{Entire Population} & \textemdash{} & \textbf{975.82} &
\makecell{\textbf{0.662}\\ \footnotesize\textbf{[0.658;0.665]}} &
\makecell{\textbf{4.402}\\ \footnotesize\textbf{[4.392; 4.413]}} &
\makecell{\textbf{3.336}\\ \footnotesize\textbf{[3.327; 3.345]}} &
\makecell{\textbf{91.416}\\ \footnotesize\textbf{[91.330; 91.503]}} &
\textbf{3.509} &
\makecell{\textbf{0.297}\\ \footnotesize\textbf{[0.293; 0.301]}} &
\textbf{0.041} \\
\bottomrule
\end{tabular}
\end{adjustbox}
\end{subtable}

\vspace{0.8em}

\begin{subtable}{\textwidth}
\centering
\caption{\textbf{Eplet level based optimization and compatibility}}
\label{tab:epleteplet_mod_sub}
\begin{adjustbox}{max width=\textwidth}
\begin{tabular}{
  L{3.4cm} C{2.0cm} C{1.6cm} C{2.2cm} C{2.7cm} C{2.7cm} C{3.4cm} C{1.8cm} C{2.35cm} C{2.60cm}}
\toprule
\multicolumn{1}{c}{\textbf{Ethnicity}} &
\multicolumn{1}{c}{\textbf{Equity weight}} &
\multicolumn{1}{c}{\textbf{Arrivals}} &
\multicolumn{1}{c}{\(\bm{F(s)}\)} &
\multicolumn{3}{c}{\textbf{HLA(s)}} &
\multicolumn{1}{c}{\(\bm{W(s)}\)} &
\multicolumn{1}{c}{\(\bm{L(s)}\)} &
\multicolumn{1}{c}{\(\bm{1 - F(s) - L^{\ast}(s)}\)} \\
\cmidrule(lr){5-7}
\multicolumn{1}{c}{} & \multicolumn{1}{c}{} & \multicolumn{1}{c}{} &
\multicolumn{1}{c}{\footnotesize(Matched)} &
\multicolumn{1}{c}{\textbf{Antigen}} &
\multicolumn{1}{c}{\textbf{Allele}} &
\multicolumn{1}{c}{\textbf{Eplets}} &
\multicolumn{1}{c}{\footnotesize(Waiting Time)} &
\multicolumn{1}{c}{\footnotesize(Left Unmatched)} &
\multicolumn{1}{c}{\footnotesize(Still in KEP)} \\
\midrule
Caucasian & 1.000 & 590.92 &
\makecell{0.652\\ \ci{0.648}{0.656}} &
\makecell{\underline{4.548}\\ \underline{\ci{4.533}{4.562}}} &
\makecell{\underline{3.303}\\ \underline{\ci{3.288}{3.319}}} &
\makecell{\underline{96.234}\\ \underline{\ci{96.128}{96.341}}} &
3.446 &
\makecell{0.290\\ \ci{0.286}{0.295}} &
0.058 \\
Afroamerican & 1.001 & 169.98 &
\makecell{0.647\\ \ci{0.640}{0.654}} &
\makecell{\underline{3.878}\\ \underline{\ci{3.849}{3.907}}} &
\makecell{\underline{2.363}\\ \underline{\ci{2.335}{2.391}}} &
\makecell{\underline{93.576}\\ \underline{\ci{93.417}{93.735}}} &
3.574 &
\makecell{0.296\\ \ci{0.289}{0.303}} &
0.057 \\
Latin & 1.002 & 146.85 &
\makecell{\underline{0.643}\\ \underline{\ci{0.636}{0.651}}} &
\makecell{\underline{4.103}\\ \underline{\ci{4.070}{4.135}}} &
\makecell{\underline{2.563}\\ \underline{\ci{2.532}{2.594}}} &
\makecell{\underline{94.280}\\ \underline{\ci{94.080}{94.479}}} &
3.742 &
\makecell{\underline{0.301}\\ \underline{\ci{0.293}{0.310}}} &
0.055 \\
Asian & 1.004 & 55.46 &
\makecell{\underline{0.656}\\ \underline{\ci{0.645}{0.667}}} &
\makecell{\underline{3.834}\\ \underline{\ci{3.793}{3.875}}} &
\makecell{\underline{1.899}\\ \underline{\ci{1.863}{1.935}}} &
\makecell{\underline{91.593}\\ \underline{\ci{91.321}{91.865}}} &
3.532 &
\makecell{0.292\\ \ci{0.281}{0.304}} &
0.051 \\
American ind.\ / Alaska nat. & 1.000 & 7.85 &
\makecell{\underline{0.682}\\ \underline{\ci{0.658}{0.706}}} &
\makecell{\underline{4.447}\\ \underline{\ci{4.334}{4.561}}} &
\makecell{\underline{2.982}\\ \underline{\ci{2.859}{3.104}}} &
\makecell{\underline{92.190}\\ \underline{\ci{91.488}{92.893}}} &
2.896 &
\makecell{\underline{0.272}\\ \underline{\ci{0.245}{0.299}}} &
0.046 \\
Native Hawaiian\ / Pacific islander & 1.000 & 5.92 &
\makecell{\underline{0.420}\\ \underline{\ci{0.387}{0.453}}} &
\makecell{\underline{3.523}\\ \underline{\ci{3.339}{3.706}}} &
\makecell{\underline{2.102}\\ \underline{\ci{1.941}{2.262}}} &
\makecell{\underline{95.046}\\ \underline{\ci{93.847}{96.245}}} &
4.488 &
\makecell{\underline{0.468}\\ \underline{\ci{0.433}{0.503}}} &
0.113 \\
\midrule
\textbf{Entire Population} & \textemdash{} & \textbf{976.98} &
\makecell{\textbf{0.649}\\ \footnotesize\textbf{[0.645;0.652]}} &
\makecell{\textbf{4.320}\\ \footnotesize\textbf{[4.309; 4.330]}} &
\makecell{\textbf{2.942}\\ \footnotesize\textbf{[2.930; 2.954]}} &
\makecell{\textbf{95.178}\\ \footnotesize\textbf{[95.107; 95.249]}} &
\textbf{3.520} &
\makecell{\textbf{0.294}\\ \footnotesize\textbf{[0.290; 0.298]}} &
\textbf{0.057} \\
\bottomrule
\end{tabular}
\end{adjustbox}
\end{subtable}

\end{table}

Very small adjustments in weights suffice to eliminate the differences in transplant probability $F(s)$ among the subpopulations and without negatively impact the overall transplant probability. In other words, there is no effectiveness equity trade-off for transplant probability. However, matching more members of subpopulations that originally experienced lower transplant probabilities leads to lower average HLA compatibility scores for these subpopulations. Matching more members of these subpopulations appears to require matching these recipients with donors that are less compatible.

\section{Discussion and Conclusions} \label{sec:discussion}

Evidence suggests that replacing antigen based HLA compatibility paradigms by allele based and eplet based paradigms in organ allocation decisions can improve transplant outcomes. The diffusion of these new paradigms and first adoptions in deceased donor programs and in KEPs has, however, raised concerns regarding negative impacts on equity, in particular across different ethnic subpopulations. Despite the lively debate on this matter in the scientific literature, the paucity of data forms a barrier to building an evidence base until the uptake of the new paradigms in KEPs reaches a larger scale. The intrinsic complexity of the compatibility paradigms and the effectiveness and equity questions raised form a second barrier. To overcome these barriers and investigate the effects on the effectiveness and equity of KEPs, we have constructed a data set from real-life data, developed a simulation-optimization framework to measure effectiveness and equity for each of the paradigms, and developed an equity weighting based approach to diminish inequities.

The data set is constructed using actual donor and recipient data but not necessarily representative of any existing KEP. Hence, the external validity of the results for any existing KEP may be limited. However, with around thousand pairs and a ten year period, the results provide realistic insight in how the effectiveness and equity of KEPs may vary with HLA paradigms and whether the resulting inequities can be addressed.

Despite the significant correlation in HLA scores among paradigms as depicted in Figure \ref{fig:paragimcorrelations}, the effectiveness as expressed by the HLA match varies considerably with the paradigm considered in KEP optimization. When supporting the view that allele, respectively eplet level HLA matches are more accurate than antigen level matches, adopting the corresponding scores in the optimization yields significant and substantial improvements.  

Regarding the second research question, the results reveal significant and substantial differences in transplant probability, the probability of leaving the KEP without a transplant, waiting time, and HLA match among ethnic subpopulations. Moreover, the results confirm the hypothesis that these differences grow when advancing from antigen to allele and eplet level paradigms.

The results obtained for the equity weighted objective functions show that differences in access and HLA match can be eliminated without negatively impacting KEP effectiveness in terms of overall transplant probability. For the case study at hand, however, this comes at the price of lower HLA compatibilities for subpopulations for which transplant probabilities increase. 

While the weights obtained in our study have no validity beyond the KEP constructed for the purpose of this research, it is interesting that the weight increments are in the order of 0.001, suggesting that equity optimization in KEPs is a delicate matter. The magnitude of the weight increment are closely related with the normalized HLA compatibility scores included in the objective function via \ref{eq:equityweighting}, which are divided by the population size of approximately one thousand. Choosing weights in the order of magnitude of the contributions to the objective function of the (second) HLA compatibility component facilitates a trade-off with equity of transplant probability (the first component). This finding calls for future multidisciplinary research into the trade-off between transplant equity and HLA compatibility.

Another future research direction is to generate stronger evidence as more allele level and eplet level data from actual KEPs become available. Actual data may additionally provide the opportunity to replace the HLA compatibility in the objective by outcomes such as LYFT and QALYs gained. As an alternative, current instance generators, which mostly rely on the somewhat indirect and abstract notion of PRA level \citep{saidman_increasing_2006, delorme2022improved}, can be extended to include HLA distributions in relevant patient and donor populations. 

Future research may also advance by combining the objective functions considered in this research with the more complex objective functions adopted in practice. Moreover, it is valuable to extend the analysis to including altruistic donors, longer cycles, and compatible pairs. \\

{\bf Acknowledgements: } The first, third and fourth author have received funding through ANID Fondecyt Regular grant 1230361.

\bibliographystyle{apalike}
\bibliography{EJORsubmission}

\end{document}